\documentclass[a4paper,reqno]{amsart}
\usepackage{amssymb}
\usepackage{latexsym}
\usepackage{amsmath}
\usepackage{euscript}
\usepackage{graphics,color}
\usepackage[all]{xy}
\usepackage[margin=3cm]{geometry}
\usepackage{MnSymbol} 

\newcommand{\be}{\begin{equation}}
\newcommand{\ee}{\end{equation}}
\newcommand{\ba}{\begin{eqnarray}}
\newcommand{\ea}{\end{eqnarray}}
\newcommand{\baa}{\begin{eqnarray*}}
\newcommand{\eaa}{\end{eqnarray*}}
\newcommand{\bb}{\begin{thebibliography}}
\newcommand{\eb}{\end{thebibliography}}

\newcommand{\bi}[1]{\bibitem{#1}}
\newcommand{\lab}[1]{\label{#1}}
\newcommand{\re}[1]{(\ref{#1})}

\newcounter{my}
\newcommand{\he}%
   {\stepcounter{equation}\setcounter{my}%
   {\value{equation}}\setcounter{equation}0%
   }%
\newcommand{\she}%
   {\setcounter{equation}{\value{my}}%
    }%

\renewcommand\t{\tilde}

\newcommand{\olsi}[1]{\,\overline{\!{#1}}} 

\newtheorem{pr}{Proposition}
\newtheorem{cor}{Corollary}

\theoremstyle{definition}

\numberwithin{equation}{section}

\begin{document}

\title[Ramanujan OPUC]{Ramanujan's trigonometric sums and para-orthogonal polynomials on the unit circle}



\author{Alexei Zhedanov}

\address{School of Mathematics, Renmin University of China, Beijing 100872, China}

\vspace*{5mm}

\begin{abstract}
Ramanujan's trigonometric sum $c_q(n)$ can be interpreted as a set of trigonometric moments of a finite measure concentrated at primitive $q$-th roots of unity with equal masses. This gives rise to sets of corresponding polynomials orthogonal on the unit circle. We present explicit expressions of these polynomials for special values of $q$, e.g. when $q=p$ or $q=2p$ or $q=p^k$, where $p$ is a prime number.    We generalize this procedure taking the Kronecker polynomial instead of cyclotomic one. In this case the moments are expressed as finite sums of $c_q(n)$ with different $q$. 
   \end{abstract}

\maketitle

\section{Introduction}
\setcounter{equation}{0}
Ramanujan's trigonometric sum $c_M(n)$ is defined as
\be
c_M(n) = \sum_{ (s,M)=1} \exp\left(\frac{2 \pi i s n}{M} \right), \lab{c_N_def} \ee
where $n$ is an integer, $M$ is a positive integer and summation in \re{c_N_def} is performed over all $s$ coprime with $M$ and smaller than $M$.

Equivalently, Ramanujan's sum can be presented as follows. Let $\zeta_1, \zeta_2, \dots, \zeta_{\varphi(M)}$ be a set of {\it primitive} roots of unity of order $M$ and  $\varphi(n)$ is Euler's totient function, i.e. the number of integers from $1$ to $n-1$ coprime with $n$.

Recall that the roots of unity $\zeta$ of order $M$ is called the {\it primitive} if $\zeta^M=1$ but $\zeta^k \ne 1$ for $k=1,2,\dots, M-1$. Any primitive root of unity is the generator of the cyclic group $\mathbb{Z}_M$   of {\it all} $M$ roots of unity of order $M$.

For example, when $M=10$ then $\varphi(M)=4$ and the primitive roots of unity are
\be
\zeta_1=\exp\left(\frac{\pi i }{5} \right), \: \zeta_2 = \exp\left(\frac{3 \pi i }{5} \right), \: \zeta_3 =\exp\left(\frac{7 \pi i }{5} \right), \: \zeta_4 = \exp\left(\frac{9 \pi i }{5} \right)
\lab{prim_10} \ee

Ramanujan's sum can then be presented as
\be
c_M(n) = \sum_{s=1}^{\varphi(M)} \zeta_s^n
\lab{c_N_prim} \ee

The cyclotomic polynomials $C_M(z)$ are closely connected with Ramanujan sums. Recall, that for any positive integer $M$ the cyclotomic polynomial is defined as the minimal polynomial of primitive roots of $M$, i.e.
\be
C_M(z) = (z-\zeta_1) (z-\zeta_2) \dots (z-\zeta_{\varphi(M)})
\lab{C_M} \ee   
One obvious relation between cyclotomic polynomials and Ramanujan's sums follows from the Vieta formulas. Indeed, put $L=\varphi(M)$ and expand the cyclotomic as a finite sum in $z^n$:
\be
C_M(z) = z^{L} + \kappa_1 z^{L-1} + \kappa_2 z^{L-2} + \dots +  \kappa_2 z^2 + \kappa_1 z + 1.
\lab{C_M_expand} \ee
(Notice the symmetricity of the expansion coefficients which follows fro the palindromic property $C_M(z) = z^L C_M(z^{-1})$ \cite{HW}). Then 
\be
\kappa_1 = - c_M(1)=-\mu(M),
\lab{kappa_1} \ee
where $\mu(n)$ is the M\"obius function \cite{HW}. Similarly
\be
\kappa_2 = \frac{1}{2} \left(c_M^2(1) - c_M(2) \right)
\lab{kappa_2} \ee
and so on.

There is another interpretation of Ramanujan's sum which will be the main subject of our paper. Namely, we interpret $c_M(n)$ as a set of {\it trigonometric moments} $\sigma_n$ with respect to some discrete measure $d \mu(z)$ on the unit circle $|z|=1$.

Recall that the given positive measure $d \mu(\theta)$ on the unit circle one can define the trigonometric moments 
\be
\sigma_n = \frac{1}{2 \pi} \, \int_{0}^{2 \pi} e^{i \theta n} d \mu(\theta) , \quad n=0, \pm 1, \pm 2, \dots  \lab{g_n_def} \ee

For example, for the simple Lebesgue measure $d \mu(\theta) = d \theta$ all moments vanish apart from $\sigma_0$
\be
\sigma_0= 1, \quad \sigma_n =0, \; n=\pm 1, \pm 2, \dots \lab{Leb_g} \ee

Using the trigonometric moments $\sigma_n$, one can construct the Toeplitz determinants \be \Delta_n= \left |
\begin{array}{cccc} \sigma_0 & \sigma_{1} & \dots &
\sigma_{n-1}\\ \sigma_{-1}& \sigma_0 & \dots & \sigma_{n-2}\\ \dots & \dots & \dots & \dots\\
\sigma_{1-n} & \sigma_{2-n} & \dots & \sigma_0 \end{array} \right
| , \quad n=1,2,\dots \lab{Delta} \ee 
Positivity of the measure $d \mu(\theta)$ leads to the obvious property
\be
\sigma_{-n} = \bar{\sigma}_n,
\lab{sym_sigma} \ee
where $\bar{\sigma}_n$ means complex conjugate of ${\sigma}_n$.

It is well known that positivity of the measure $d \mu(\theta)$ is equivalent to positivity of  Toeplitz determinants \cite{Simon}
\be
\Delta_n >0, \quad n=1,2,3, \dots \lab{pos_Delta} \ee

Given the measure $\mu(\theta)$ one can introduce the polynomials \cite{Simon}
\ba \nonumber
&&\Phi_n(z)=(\Delta_n)^{-1} \left |
\begin{array}{cccc} \sigma_0 & \sigma_1 & \dots & \sigma_n \\ \sigma_{-1} & \sigma_0 &
\dots & \sigma_{n-1} \\ \dots & \dots & \dots & \dots\\
\sigma_{1-n}& \sigma_{2-n}& \dots & \sigma_1\\ 1& z & \dots & z^n
\end{array} \right |, \lab{deterPhi} \ea
It is seen that the polynomial $\Phi_n(z)$ is a $n$-degree monic polynomial
\be
\Phi_n(z) = z^n + O(z^{n-1}) \lab{Phi_monic} \ee
By construction, the polynomials $\Phi_n(z)$ possess the orthogonality property
\be
 \int_{0}^{2 \pi} \Phi_n(e^{i \theta}) e^{-i j \theta } d \mu(\theta) = 0, \; j=0,1,\dots, n-1
 \lab{ort_1} \ee
which is equivalent to \cite{Simon}
\be
\int_{0}^{2 \pi} \Phi_n(e^{i \theta}) \olsi{\Phi_m(e^{i \theta})} d \mu(\theta) = \int_{0}^{2\pi} \Phi_n(e^{i \theta}) \olsi{\Phi}_m(e^{-i \theta}) d \mu(\theta)= h_n \, \delta_{nm}, \lab{ort_2} \ee
where
\be
h_n = \frac{\Delta_{n+1}}{\Delta_n}>0. \lab{h_Delta} \ee
(Notation $\olsi{\Phi}_m(z)$ means that we take only conjugate of the coefficients of the polynomial $\Phi_n(z)$ but leave the argument $z$ unchanged).
 
The polynomials $\Phi_n(z)$ are called the polynomials orthogonal on the unit circle, or OPUC for brevity (we adopt the abbreviation of \cite{Simon}).

They satisfy the fundamental Szeg\H{o} recurrence relation \cite{Simon}
\be
\Phi_{n+1}(z) = z \Phi_n(z) - \bar a_n \Phi_n^*(z), \lab{Sz_rec} \ee 
where
\be
\Phi_n^*(z) = z^n \olsi{\Phi}_n(z^{-1}).
\lab{Phi*} \ee
The parameters $a_n$ (sometimes called the Verblusnky parameters \cite{Simon}) satisfy the condition
\be
|a_n| <1 . \lab{a<1} \ee
One can show that the recurrence relation \re{Sz_rec} together with the condition \re{a<1} is necessary and sufficient for polynomials $\Phi_n(z)$ to be OPUC with respect to a positive measure $d \mu(\theta)$ \cite{Simon}. The normalization constants $h_n$ have the following expression in terms of $a_n$ \cite{Simon}:
\be
h_0=1, \; h_n = \left(1-|a_0|^2 \right) \left(1-|a_1|^2 \right) \dots \left(1-|a_{n-1}|^2 \right), \quad n=1,2,\dots
\lab{h_a} \ee

There is an important special case when the number of OPUC is {\it finite}. This happens if $|a_i|<1$ for $i=0,1,\dots, N-1$ but $|a_N|=1$. Then the polynomial $\Phi_{N+1}(z)$ has $N+1$ simple roots $z_s, \: s=1,2,\dots, N+1$ on the unit circle and one has the orthogonality relation
\be
\sum_{s=1}^{N+1} \Phi_n(z_s) {\bar \Phi}_m({z_s}^{-1}) w_s = h_n \delta_{nm}, \lab{fin_ort} \ee
where the positive weights are \cite{Ger}, \cite{Simon}
\be
w_s = \frac{h_N}{\Phi_{N+1}'(z_s) {\bar \Phi}_N(z^{-1}_s)} = \frac{h_N}{\olsi{\Phi_{N+1}'(z_s)}  \Phi_N(z_s)}.
\lab{w_s_OPUC} \ee
(the last equality in \re{w_s_OPUC} holds because the product $\Phi_{N+1}'(z_s) {\bar \Phi}_N(z^{-1}_s)$ is real). 
Such finite system of orthogonal polynomials is called the para-orthogonal polynomials on the unit circle (POPUC) \cite{Simon}. 

In this paper we propose new explicit examples of finite systems POPUC with equal concentrated masses located at primitive roots of unity. More exactly, for a given positive integer $M$ we put $N=\varphi(M)$ and define the system $\Phi_0(z)=1, \Phi_1(z), \dots, \Phi_N(z), \Phi_{N+1}(z)$ of POPUC by identifying their trigonometric moments $\sigma_n$ with the Ramanujan sum 
\be
\sigma_n = c_M(n), \quad n=0,1,2,\dots, N+1
\lab{sigma_c} \ee
The first (obvious) proposition is about the measure corresponding to the choice \re{sigma_c}
\begin{pr}
The moments \re{sigma_c} correspond to the finite discrete measure on the unit circle located at primitive roots of unity $\zeta_1, \zeta_2, \dots, \zeta_{N+1}$ with equal concentrated masses:
\be
w_s = (N+1)^{-1}.
\lab{w-equal} \ee
\end{pr}
Proof of this proposition is trivial. Indeed, by definition of the moments we have
\be 
\sigma_n = \sum_{s=1}^{N+1} w_s \zeta_s^n = (N+1)^{-1} \sum_{s=1}^{N+1} \zeta_s^n
\lab{sigma_dis} \ee 
which coincide (up to a trivial factor) with definition \re{c_N_prim} of the Ramanujan's sum.

It is natural to call the corresponding POPUC $\Phi_n(z)$ the {\it Ramanujan para-orthogonal polynomials}. 

The main problem will be to find explicitly the POPUC $\Phi_n(z)$ and corresponding Verblunsky parameters $a_n$. This will be done in the next sections. Moreover, we prove an important result concerning relations between "Sturmian" POPUC and Ranaujan POPUC. Namely, we show that these two systems are mirror-dual one with respect to another. This is the main result of the Section 2. In Section 3 and 4 we present several explicit examples of such dual systems. They are connected either with the cyclotomic or with the Kronecker polynomials.

\section{Reconstruction of polynomials orthogonal on the unit circle}
\setcounter{equation}{0}
In order to reconstruct the polynomials $\Phi_n(z)$, we use the following proposition which plays the main role in whole theory
\begin{pr}
Assume that $\Phi_0(z), \Phi_1(z), \dots, \Phi_N(z), \Phi_{N+1}(z)$ is a set of POPUC corresponding to the Verblunsky parameters $a_0,a_1, \dots, a_{N-1}, a_N$, where $|a_k|<1, \: k=0,1,\dots, N-1$ and $|a_N|=1$. Assume that these polynomials are orthogonal on the unit circle at some points $z_k$ with concentrated masses $w_s$, i.e.
\be
\sum_{s=1}^{N+1} \Phi_n(z_s) \olsi {\Phi_m(z_s)}w_s =(N+1) h_n \delta_{nm}
\lab{gen_eqw} \ee

Let $\t \Phi_n(z)\: n=0,1,\dots, N+1$ be another system of POPUC (called the "mirror-dual"). This system is characterized by the property that their Verblunsky parameters $\t a_n$ are "mirror-dual" with respect to $a_n$. This means that \cite{Mar1}
\be
\t a_n = - a_N {\bar a}_{N-n-1}, \quad n=0,1,\dots, N,
\lab{ta_def} \ee
where it is assumed that $a_{-1}=-1$. (Note that  the coefficient $a_{-1}$ does not play any role in expressions of the polynomials $\Phi_n(z)$, so we can choose it arbitrarily, however, it is naturally to put $a_{-1} = -\bar \Phi_0(z) =-1$).

Then the polynomials $\t \Phi_n(z)$ are orthogonal on the {\it same} set of points $z_1,z_2, \dots, z_{N+1}$  
\be
\sum_{s=1}^{N+1} \t \Phi_n(z_s) \olsi {\t \Phi_m(z_s)} \t w_s = \t h_n \delta_{nm}
\lab{dual_ort} \ee 
with the concentrated positive masses $\t w_s>0$ determined uniquely from the relation
\be
\t w_s w_s = \frac{h_{N}}{|\Phi'_{N+1}(z_s)|^2},  
\lab{tww} \ee
where $h_n$ are given by \re{h_a}.

Conversely, assume that the weights $w_s$ and $\t w_s$ of two different systems $\Phi_n(z), \, \t \Phi_n(z)$ of POPUC have the same spectral points $z_s$ on the unit circle  and their concentrated masses  are related as in \re{tww}. Then these systems are mirror-dual and their Verblunsky parameters are related as in \re{ta_def}. 

\end{pr}

{\it Remark 1}. In fact, condition \re{tww} can be replaced with a more weak condition 
\be
\t w_s w_s = \frac{1}{|\Phi'_{N+1}(z_s)|^2},  
\lab{tww1} \ee
because the weights $w_s, \t w_s$ are determined up to an arbitrary normalization factor. Note also that 
\be
\t h_N = h_N ,
\lab{thh} \ee
which follows trivially from \re{h_a}.

\vspace{4mm}

{\it Remark 2}. Condition \re{tww} (or \re{tww1}) is an exact "circle" analog of the corresponding relation between two finite sets of the mirror-dual orthogonal polynomials on the real line, see \cite{Borodin} and \cite{VZ_dual} for details.

\vspace{4mm}

{\it Proof}. In fact, the above proposition follows easily from the Theorem {\bf 7} in \cite{Mar1}. Namely, it was showed in \cite{Mar1} that if the Verblunsky parameters satisfy the mirror-dual relation \re{ta_def} then their "final" characteristic polynomials coincide
\be
\t \Phi_{N+1}(z) = \Phi_{N+1}(z) \lab{t_N+1} \ee
(and hence the spectral points $z_s$ on the unit circle are the same). Moreover, there is the formula \cite{Mar1}
\be
\frac{\Phi_N(z)}{\Phi_{N+1}(z)} = \sum_{s=1}^{N+1} \frac{\t w_s}{z-z_s}, \lab{t_frac} \ee
where the weights are appropriately normalized
\be
\sum_{s=1}^{N+1} \t w_s =1.
\lab{norm_tw} \ee
From \re{t_frac} it follows that
\be
\t w_s = \frac{\Phi_N(z_s)}{\Phi_{N+1}'(z_s)}.
\lab{tws_Phi} \ee
Now, multiplying \re{w_s_OPUC} and \re{tws_Phi} we arrive at \re{tww}.

Conversely, assume that the (positive) weights $\t w_s$ are determined by \re{tww1} on the same spectral points $z_s$. Then the whole system of monic POPUC $\t \Phi_n(z)$ is determined uniquely. Hence the Verblunsky parameters $\t a_n$ should satisfy the mirror-dual condition \re{ta_def}.

There is an easy consequence of this proposition which is important for our further construction of Ramanujan POPUC.

\vspace{4mm}

\begin{cor} 
The POPUC $\Phi_n(z)$ with equal masses $w_s=(N+1)^{-1}$ on $N+1$ arbitrary distinct points $z_s, \: s=1,2,\dots, N+1$ on the unit circle are mirror-dual with respect to so-called Sturmian POPUC $\t \Phi_n(z)$. The latter are characterized by the condition \cite{Zhe_K}
\be
\t \Phi_N(z) = (N+1)^{-1} \t \Phi_{N+1}'(z).
\lab{Sturm_cond} \ee
\end{cor}
The proof is easy. Indeed, if all masses $w_s$ are the same $w_s=(N+1)^{-1}$, then from \re{tww} we have 
\be
\frac{\t \Phi_N(z)}{\t \Phi_{N+1}(z)} = \sum_{s=1}^{N+1} \frac{w_s}{z-z_s} = (N+1)^{-1} \sum_{s=1}^{N+1} \frac{1}{z-z_s}
\lab{pft} \ee
this leads to 
\be
\t \Phi_N(z) = (N+1)^{-1} \t \Phi_{N+1}'(z) = (N+1)^{-1} \Phi_{N+1}'(z) \lab{Sturm_cond2} \ee
which gives \re{Sturm_cond}.

On the other hand the Sturmian POPUC $\t \Phi_n(z)$ are completely characterized by condition \cite{Zhe_K}
\be
\t w_s = \frac{const}{|\Phi'_{N+1}(z_s)|^2}.
\lab{tw_Sturm} \ee 
Hence the mirror-dual POPUC have the equal concentrated masses.

This Corollary allows to find the Verblunsky parameters $a_n$ of the Ramanujan POPUC from the already known parameters $\t a_n$ for the Sturmian cyclotomic POPUC considered in \cite{Zhe_K}. 

{\it Remark}. For the case of polynomials on the real line this Corollary was exploited in \cite{Zhe_S}.

\section{Simplest cases of Ramanujan cyclotomic POPUC}
\setcounter{equation}{0}
It is natural to take the polynomial $\Phi_{N+1}(z)$ as a cyclotomic polynomial $C_M(z)$. Indeed, in this case all the spectral points $z_k$ coincide with the primitive roots of unity. Hence the moments are $\sigma_n=(N+1)^{-1}c_M(n)$.  The simplest case of the cyclotomic polynomial $\Phi_{N+1}(z)=C_{p}(z)$ where $p$ is an arbitrary odd prime. In this case $N=p-2$ and
\be
\Phi_{N+1} = \frac{z^p-1}{z-1} = z^{p-1} + z^{p-2} + \dots + z+1.
\lab{Phi_N+1_p} \ee 

The moments are
\be
\sigma_0=1, \quad \sigma_n = -(p-1)^{-1}, \; n=1,2,\dots, p-1
\lab{sigma_prime} \ee
With these moments one can easily calculate all the determinants \re{Delta} and \re{deterPhi} to get
\be
\Delta_n = \frac{p^{n-1}(p-n)}{(p-1)^n}
\lab{Delta_p} \ee
and
\be
\Phi_n(z) =  z^n + (N-n+2)^{-1} \left(z^{n-1} + z^{n-2}+ \dots + z+1 \right), \quad n=0,1,2, \dots, N, N+1
\lab{Phi_prime_p} \ee
The Verblunsky parameters are
\be
a_n  =-\Phi_{n+1}(0) = -\frac{1}{N-n+1}, \quad n=0,1,\dots, N .
\lab{a_prime_p} \ee

The Sturmian cyclotomic POPUC were identified in \cite{Zhe_K} with the well known "single-moment" polynomials. They have the explicit expression 
\be
\t \Phi_n(z)= \frac{1}{n+1} \sum_{k=0}^n (k+1)z^k, \quad n=0,1, \dots, N
\lab{smomP} \ee
and corresponding Verblusnky parameters are
\be
\t a_n= -\frac{1}{n+2}, \; n=0,1,2,\dots, N-1, \; a_N =-1. \lab{a_single} \ee
As expected, the Verblunsky coefficients \re{a_prime_p} and \re{a_single} are related with the mirror-dual correspondence \re{ta_def}. 

The case of the cyclotomic polynomials with $M=2p$ i.e. with the condition $\Phi_{N+1}(z)=C_{2p}(z)$ is also easy. In this case the moments $\sigma_n$ differ from the previous example \re{sigma_prime} only by sign changing behavior:
\be
\sigma_0=1, \quad \sigma_n = (-1)^{n+1} (p-1)^{-1}, \; n=1,2,\dots, p-1
\lab{sigma_prime_a} \ee
Hence the Verblunsky coefficients
\be
a_n = (-1)^n \frac{1}{N-n+1}, \quad n=0,1,\dots, N \lab{a_2p} \ee

Similarly, one can consider the case when $M = 2^j p^k$, where $p$ is an odd prime number and $k,k$ are arbitrary positive integers. A showed in \cite{Zhe_K} this leads to a set of {\it sieved} (see, e.g. \cite{Ismail} on definition of sieved OPUC) POPUC. 

More complicated is the case $M= p_1 p_2 \dots p_m$ when $M$ is a product of distinct odd prime numbers $p_1<p_2<
\dots <p_m$. Already the simplest - so called {\it binary} - case, when $M=pq$ with only two primes $p<q$, is very nontrivail. We were able only to formulate a conjecture concerningg behavior  for the "head" and "tail" of the sequence of Verblunsky parameters $a_n$ \cite{Zhe_K}. Clearly, this conjecture remains valid for the Ramanujan POPUC with interchanging "head" and "tail".

\section{Ramanujan Kronecker POPUC}
\setcounter{equation}{0}
The Kronecker polynomial $K(z) = z^M + e_{M-1}z^{M-1} + e_{M-2}z^{M-2} + \dots + e_1 z + e_0$  is a monic polynomial with integer coefficients $e_k$ such that all its roots $z_k$ lie in the closed unit disc $|z_k|<1$ \cite{Dam}. It is easy to show if $a_0 \ne 0$ then any Kronecker polynomial has all roots belonging to the unit circle $|z_k|=1$ and that in fact any Kronecker polynomial (with nonzero roots) can be presented as a finite product of cyclotomic polynomials
\be
K(z) = C_{m_1}^{j_1}(z) C_{m_2}^{j_2}(z) \dots C_{m_k}^{j_k}(z) 
\lab{K_C} \ee    
with some positive integers $m_i,\,  j_i$. 

If the polynomial $K(z)$ is taken as the polynomial $K(z) =\Phi_{N+1}(z)$ belonging to a set of POPUC $\Phi_0(z), \Phi_1(z), \dots, \Phi_N(z), \Phi_{N+1}(z)$, then all roots of $\Phi_{N+1}(z)$ should be simple \cite{Mar1}. This leads to the presentation of all such admissible polynomials:
\be
K(z) =  C_{m_1}(z) C_{m_2}(z) \dots C_{m_k}(z), 
\lab{K_C_adm} \ee
where all $m_i$ are distinct positive inters.

Indeed, if all $m_i$ are distinct, then the set of roots of $K(z)$ is a disjoint union of primitive roots of unity.

We can define the {\it Ramanujan Kronecker polynomials} as a set of POPUC  $\Phi_0(z), \Phi_1(z), \dots, \Phi_N(z), \Phi_{N+1}(z)$ such that 
$\Phi_{N+1}(z) = K(z)$ with $K(z)$ given by \re{K_C_adm} and with {\it equal concentrated masses} at all roots of $\Phi_{N+1}(z)$.

Let us denote $z_1, z_2, \dots z_M$ all the roots of the polynomial $K(z)$. Then the orthogonality relation for the Ramanujan Kronecker POPUC reads
\be
\sum_{s=1}^{M} \Phi_n(z_s) \Phi_m(z_s^{-1}) = M h_n \delta_{mn}, \quad m,n=0,1,\dots , N
\lab{ort_Kron} \ee

We have the 
\begin{pr}
The moments $\sigma_n$ of the Ramanujan Kronecker polynomials have the presentation
\be
\sigma_n = M^{-1} \left(c_{m_1}(n) + c_{m_2}(n) + \dots + c_{m_k}(n) \right)
\lab{sigma_Ram} \ee
\end{pr}

{\it Proof}. We have that the concentrated masses at all roots of $K(z)$ are the same. Hecne, by definition of  moments 
\be
\sigma_n = \sum_{s=1}^M \kappa \,  z_s^n,
\lab{sigma_Ram_def} \ee
where $\kappa$ is an appropriate normalization factor. Obviously, $\kappa = M^{-1}$. And we can perform summation in \re{sigma_Ram_def} step-by-step: first performing summation over all the roots of $C_{m_1}$, then of $C_{m_2}(z)$ and so on.   Each such summation involves moments of cyclotomic polynomials and hence we arrive at \re{sigma_Ram}.

We thus have a simple alternative definition of Ramanujan Kronecker POPUC: these are POPUC with the moments $\sigma_n$ given by a finite sum of different Ramanujan's sums $c_m(n)$ \re{sigma_Ram}.

We consider here only one special case of the Ramanujan Kronecker polynomials which is generated by the polynomial $K(z) = A_{2p}(z)$, where $A_{M}(z)$ is so-called anti-cyclotomic (or inverse cyclotomic) polynomial defined in \cite{Moree} as Kronecker polynomials with roots being complementary to primitive roots of unity of order $M$. In other words, the polynomial $A_M(z)$ has the expression
\be
A_{M}(z) = \frac{z^M-1}{C_M(z)} \lab{A_M_def} \ee

In \cite{Zhe_K} we have considered the Sturmian Kronecker POPUC corresponding to the anti-cyclotomic polynomial with $p$ an arbitrary odd prime number. In this case $N=p$ and 
\be
\Phi_{N+1}(z) = A_{2p}(z) = z^{p+1} + z^{p} -z-1
\lab{A_{2p}_exp} \ee 
The Sturmian POPUC are defined by the condition
\be
\Phi_N(z) = \frac{A_{2p}'(z)}{p+1}=
    \frac{z^p +pz^{p-1}-1}{p+1}
\lab{Phi_N_anti_C} \ee
Then further polynomials $\Phi_k(z), \: k=N-1,N-2, \dots, 0$ are determined uniquely by the inverse Szeg\H{o} recurrence relation \cite{Zhe_K}.

The results are \cite{Zhe_K}
\be
\t \Phi_0(z)=1, \: \t \Phi_1(z) = z+ 1-p^{-1}, \quad \t \Phi_n(z) = {z}^{n}+{\frac { \left( 2\,p-n \right) {z}^{n-1}}{2\,p-n+1}}+{\frac {
 \left( -1 \right) ^{n}}{2\,p-n+1}}, \; n=2,3,\dots, p.
\lab{Phi_anti} \ee 
with the Verblunsky parameters
\be 
\t a_0 = \frac{1-p}{p}, \quad \t a_n = \frac{(-1)^{n}}{2p-n}, \; n=1, 2, \dots, p-1, \quad a_{p}=1.
\lab{a_n_A} \ee

Using mirror-reflection rule we arrive at the following expression for the Verblunsky parameters of the Ramanujan anti-cyclotomic polynomials 
\be
a_n = - \t a_{p-n-1} ,\quad  n=0,\dots, p
\lab{a_n_RA} \ee
It is easily verified that corresponding Ramanujan anti-cyclotomic POPUC have the expression
\be
\Phi_n(z) = z^n + \frac{1}{n+p} \, \frac{z^n -(-1)^n}{z+1}, \quad n=0,1,\dots, p. \lab{Phi_n_AR} \ee

\vspace{5mm}

\section{Conclusions}
\setcounter{equation}{0}
The main result of \cite{Zhe_K} and of this paper is a possibility to construct a wide family of finite para-orthogonal polynomials starting with prescribed cyclotomic or with more general Kronecker polynomials.  These polynomials are taken as the final members $\Phi_{N+1}(z)$ of the set $\Phi_0, \Phi_1(z), \dots, \Phi_N(z)$. The two systems presented in \cite{Zhe_K} and in the present paper differ by the algorithm of construction. In \cite{Zhe_K} we start with two polynomials: $\Phi_{N+1}(z)$ and $\Phi_N(z)$ taken as the derivative of $\Phi_{N+1}(z)$. Then all further polynomials can be reconstructed by a specific "circle Sturm" algorithm, like for polynomials on the real line.

In the present paper another algorithm is proposed: we start with the same polynomial $\Phi_{N+1}(z)$ but demand that the concentrated masses be equal at all spectral points which are roots of $\Phi_{N+1}(z)$. The main result is presented in the Corollary 1: it appears that these two systems of POPUC are mirror-dual.  Another interesting result is possibility to identify the famous Ramanujan' trigonometric sums as moments of corresponding POPUC.    

We have presented several explicit examples of the polynomials $\Phi_n(z)$. It is expected that many new explicit examples can be constructed by the above algorithms.

\vspace{6mm}

{\large \bf Acknowledgments.} 

\vspace{5mm}

 The author thanks L.Vinet for illuminating communications.  The author is gratefully holding Simons CRM Professorship. The work is funded by the National Foundation of China (Grant No.11771015).

\bb{99}

\bi{Borodin}  A. Borodin, {\it Duality of orthogonal polynomials on a finite set}, J. Statist. Phys. {\bf 109} (2002) 1109--1120.




\bi{Dam} P.A.Damianou, {\it Monic Polynomials in $Z[x]$ with Roots in the Unit Disc}, Am.Math.Month., {\bf 108} (2001), 253--257.

\bi{Ger} Ya.L. Geronimus,\quad {\it Polynomials Orthogonal on a
Circle and their Applications}, \\ Am.Math.Transl.,Ser.1 {\bf
3}(1962), 1-78.


\bi{HW} G.Hardy and E.M.Wright, An introduction to the theory of numbers (Sixth ed.), 2008, Oxford University Press.


\bi{Ismail} M.E.H.Ismail and Xin Li, {\it On sieved orthogonal polynomials IX: Orthogonality on the unit circle}, Pacific J.Math. {\bf 153}, (1992), 289--297.




\bi{Mar1} A.~Mart\'inez-Finkelshtein, B.~Simanek, and B.~Simon. {\it Poncelet's theorem, paraorthogonal polynomials
and the numerical range of compressed multiplication operators}. Adv. Math., {\bf 349}, 992--1035, (2019).

\bi{Moree} P.Moree, {\it Inverse cyclotomic polynomials}, J.Numb.Th. {\bf 129} (2009) 667--680.

\bi{Simon} B.Simon, {\it Orthogonal Polynomials On The Unit
Circle}, AMS, 2005.

\bi{VZ_dual} L.Vinet and A.Zhedanov, {\it A characterization of classical and semiclassical orthogonal
polynomials from their dual polynomials}, J.Comp.Appl.Math. {\bf 172}, (2004) 41--48.



\bi{Zhe_S} A.Zhedanov, {\it Classical Sturmian sequences}, arXiv:1904.03789.

\bi{Zhe_K} A.Zhedanov, {\it Para-orthogonal polynomials on the unit circle generated by Kronecker polynomials}, arXiv:2107.11430.

\eb

\end{document}